\documentclass[notitlepage]{amsart}

\usepackage{amsthm, amssymb, latexsym, tabls, hhline,lscape}
\input xypic
\xyoption{all} \xyoption{arc}

\newtheorem{prop}{Proposition}[section]

\newtheorem{thm}[prop]{Theorem}
\newtheorem{cor}[prop]{Corollary}
\newtheorem{lem}[prop]{Lemma}

\def\vertical{$\vphantom{\frac{\DS{A}}{\DS{A}}}$}
\def\verticala{\vphantom{\frac{\DS{A}}{\DS{A}}}}
\def\equat{\refstepcounter{prop}$$~}
\def\endequat{\leqno{\boldsymbol{(\arabic{section}.\arabic{prop})}}~$$}
\def\exemple#1{{\refstepcounter{prop}\label{#1}\noindent\it Example
\arabic{section}.\arabic{prop} - }}

\def\remarque#1{{\refstepcounter{prop}\label{#1}\noindent\it Remark
\arabic{section}.\arabic{prop} - }}

\def\sba{{\bar{s}}}

\mathchardef\inferieur="321E

\def\DC{{\mathcal{D}}}

\def\PC{{\mathcal{P}}}

\def\RM{{\mathbb{R}}}
\def\NM{{\mathbb{N}}}
\def\QM{{\mathbb{Q}}}
\def\ZM{{\mathbb{Z}}}

\def\CM{{\mathbb{C}}}

\def\UD{{\mathcal{D}}}

\def\UP{{\mathcal{P}}}
\def\itemth#1{\item[$({\mathrm{#1}})$]}
\def\e{\varepsilon}

\newcommand{\card}[1]{{\left\arrowvert #1 \right\arrowvert\,}}

\newcommand{\set}[1]{{\left\{ #1 \right\}}}  
\def\st{{\,\arrowvert\, }}  

\def\a{\alpha}

\def\g{\gamma}

\def\D{\Delta}
\def\e{\varepsilon}

\def\r{\rho}

\def\th{\theta}

\def\z{\zeta}

\def\iff{\Leftrightarrow}

\def\to{\rightarrow}

\def\lexp#1#2{\kern\scriptspace\vphantom{#2}^{#1}\kern-\scriptspace#2}
\def\finl{~$\scriptstyle{\square}$}
\def\SS{\scriptstyle}
\def\DS{\displaystyle}

\def\SG{{\mathfrak S}}

\def\Rad{{\textnormal{Rad}}}
\def\ad{{\textnormal{ad}}}
\def\Ind{\mathop{\mathrm{Ind}}\nolimits}
\def\Ker{\mathop{\mathrm{Ker}}\nolimits}
\def\rank{\mathop{\mathrm{rank}}\nolimits}

\def\edclass{{\ \sim_{A}\ }} 

\def\Irr{\mathop{\mathrm{Irr}}\nolimits}

\def\itemth#1{\item[$({\mathrm{#1}})$]}
\def\finl{~$\SS{\square}$}

\begin{document}

\title{Generalized descent algebras}
\author{C\'edric Bonnaf\'e}
\address{Universit\'e de Franche-Comt\'e,
D\'epartement de Math\'ematiques,
16 route de Gray,
25 000 Besan\c con, France}
\email{bonnafe@descartes.univ-fcomte.fr}

\author{Christophe Hohlweg$^1$}
\address{The Fields Institute,
222 College Street, Toronto, Ontario, Canada M5T 3J1}
\email{chohlweg@fields.utoronto.ca}

\date{\today}

\maketitle

\footnotetext[1]{Partially supported by Canada Research Chairs}

\begin{abstract}
If $A$ is a subset of the set of reflections of a finite Coxeter
group $W$, we define a sub-$\ZM$-module $\UD_A(W)$ of the group
algebra $\ZM W$. We provide examples where this submodule is a
subalgebra. This family of subalgebras includes strictly the
Solomon descent algebra and, if $W$ if of type $B$, the
Mantaci-Reutenauer algebra.
\end{abstract}

\markboth{}{}

\section*{Introduction}

Let $(W,S)$ be a finite Coxeter system whose length function is
denoted by $\ell: W \to \NM=\{0,1,2,\dots\}$. In 1976, Solomon
introduced a remarkable subalgebra $\Sigma W$ of the group algebra
$\ZM W$, called the {\it Solomon descent algebra} \cite{solomon}.
Let us recall its definition. If $I \subset S$, let $W_I$ denote
the {\it standard parabolic subgroup} generated by $I$.
Then
$$
X_I =\set{w\in W \st \,\forall\,s \in I,\,\ell(ws)>\ell(w)}
$$
is a set of {\it minimal length coset representatives} of $W/W_I$.
Let $x_I=\sum_{w\in X_I} w \in \ZM W$. Then $\Sigma W$ is defined
as the sub-$\ZM$-module of $\ZM W$ spanned by $(x_I)_{I \subset
S}$. Moreover, $\Sigma W$ is endowed with a $\ZM$-linear map $\th:
\Sigma W \to \ZM\Irr W$ satisfying $\th(x_I)=\Ind_{W_I}^W
1_{W_I}$. This is an algebra homomorphism.

If $W$ is the symmetric group $\SG_n$, $\theta$ becomes an
epimorphism and the pair $(\Sigma \SG_n,\th)$ provides a nice
construction of $\Irr(\SG_n)$ \cite{jollen}, which is the first
ingredient of several recent works, see for instance
\cite{Thibon,bless}. However, the morphism $\theta$ is surjective
if and only if $W$ is a product of symmetric groups.

In \cite{bonnafe-hohlweg}, we have constructed a subalgebra
$\Sigma'(W_n)$ of the group algebra $\ZM W_n$ of the Coxeter group
$W_n$ of type $B_n$: it turns out that this subalgebra contains
$\Sigma W_n$, that it is also endowed with an algebra homomorphism
$\th': \Sigma'(W_n) \to \ZM\Irr W_n$ extending $\th$. Moreover,
$\th'$ is surjective and $\QM
\otimes_\ZM \Ker \th'$ is the radical of $\QM \otimes_\ZM
\Sigma'(W_n)$. This lead to a construction of the irreducible
characters of $W_n$ following J\H ollenbeck's strategy. In fact,
$\Sigma'(W_n)$ is the {\it Mantaci-Reutenauer algebra}
\cite{mantaci-reutenauer}.

It is a natural question to ask whether this kind of construction
can be generalized to other groups.  However, the situation seems
to be much more complicated in the other types. Let us explain now
what kind of subalgebras we are looking for.

Let $T=\{wsw^{-1}~|~w \in W$ and $s \in S\}$ be the set of
reflections in $W$. Let $A$ be a fixed subset of $T$. If $I
\subset A$, we still denote
by $W_I$ the subgroup of $W$ generated by $I$ and we still set
$$
X_I =\set{w\in W \st \,\forall\,x \in W_I,\,\ell(wx) \ge \ell(w)}.
$$
Then $X_I$ is again a set of representatives for $W/W_I$. Now, let
$x_I=\sum_{w \in X_I} w \in \ZM W$. Then
$$
\Sigma_A(W)=\sum_{I \subset A} \ZM x_I
$$
is a sub-$\ZM$-module of $\ZM W$. However, it is not in general a
subalgebra of $\ZM W$.

Let us define another sub-$\ZM$-module of $\ZM W$. If $w \in W$,
let
$$
D_A(w)=\{s \in A~|~\ell(ws) < \ell(w)\} \quad\subset A
$$
be the $A$-descent set of $w$. A subset $I$ of $A$ is said to be
{\it $A$-admissible} if there exists $w \in W$ such that
$D_A(w)=I$. Let $\PC_\ad(A)$ denote the set of $A$-admissible
subsets of $A$. If $I \in \PC_\ad(A)$, we set $D_I^A=\{w \in
W~|~D_A(w)=I\}$ and $d_I^A=\sum_{w \in D_I^A} w \in \ZM W$. Now,
let
$$
\UD_A(W)=\mathop{\oplus}_{I \in \PC_\ad(A)} \ZM d_I^A.
$$
As an example, $\UD_S(W)=\Sigma_S(W)=\Sigma W$.
The main theorem of this paper is the following (here, $C(w)$
denotes the conjugacy class of $w$ in $W$):

\medskip

\noindent{\bf Theorem A.}
{\it If there exists two subsets $S_1$ and $S_2$ of $S$ such that
$A=S_1 \cup (\cup_{s \in S_2} C(s))$, then $\UD_A(W)$ is a
subalgebra of $\ZM W$.}

\medskip

This theorem contains the case of Solomon descent algebra (take
$A=S$), the Mantaci-Reutenauer algebra $\Sigma'(W_n)$ (take
$A=\{s_1,\dots,s_{n-1}\} \cup C(t)$, where
$S=\{t,s_1,\dots,s_{n-1}\}$ satisfies $C(t) \cap
\{s_1,\dots,s_{n-1}\}=\varnothing$), and
one of the subalgebras we have constructed in type $G_2$. But it
also gives a new algebra in type $F_4$ of $\ZM$-rank $300$ (take
$A=\{s_1,s_2\} \cup C(s_3)$, where $S=\{s_1,s_2,s_3,s_4\}$
satisfies $\{s_1,s_2\} \cap C(s_3)=\varnothing$). Moreover, if
$A=T$, we get that $\UD_A(W)=\ZM W$ (see Example \ref{group
algebra}). In the case of dihedral groups, we get another family
of algebras:

\medskip

\noindent{\bf Theorem B.}
{\it If $W$ is a dihedral group of order $4m$ ($m\geq 1$),
$S=\set{s,t}$ and $A=\set{s,t,sts}$ or $A=\set{t,sts}$, then
$\UD_A(W)$ is a subalgebra of $\ZM W$.}

\medskip

It must be noted that the algebras constructed in Theorems A and B
are not necessarily unitary. More precisely $1 \in \UD_A(W)$ if
and only if $S\subset A$. Moreover, if $S \subset A$, then $\Sigma
W \subset \UD_A(W)$.

The bad point in the above construction is that, in many of the
above cases, one has $\Sigma_A(W) \not= \UD_A(W)$. Also, we are
not able to construct in general a morphism of algebras $\th_A:
\UD_A(W) \to \ZM \Irr W$ extending $\th$ if $S \subset A$. It can
be proved that, in some cases, this morphism does not exist.

This paper is organized as follows. Section 1 is essentially
devoted to the proofs of Theorems A and B. One of the key step is
that the subsets $D_I^A$ are left-connected (recall that a subset
$E$ of $W$ is said to be {\it left-connected} if, for all $w$, $w'
\in E$, there exists a sequence $w=w_1$, $w_2$,\dots, $w_r=w'$
of elements of $E$ such that $w_{i+1} w_i^{-1} \in S$ for every $i
\in \{1,2,\dots,r-1\}$). In Section 2, we discuss more precisely
the case of dihedral groups.

%
%
%
%
%

\section{Descent sets}

Let $(W,S)$ be a finitely generated Coxeter system (not necessary
finite). If $s$, $s' \in S$, we denote by $m(s,s')$  the order of
$ss'\in W$. If $W$ is finite, we denote by $w_0$ its longest
element.

\subsection{Root system}
Let $V$ be an $\RM$-vector space endowed with a basis indexed by $S$ 
denoted by $\D=\set{\a_s \st s\in S}$. Let $B: V \times V \to
\RM$ be the symmetric bilinear form such that
$$
B(\a_s,\a_{s'}) = - \cos\left(\frac{\pi}{m(s,s')}  \right)
$$
for all $s$, $s' \in S$. If $s \in S$ and $v \in V$, we set
$$
s (v)=v-2B(\a_s , v)\a_s .
$$
Thus $s$ acts as the reflection in the hyperplane
orthogonal to $\alpha_s$ (for the bilinear form $B$). This extends
to an action of $W$ on $V$ as a group generated by reflections. It
stabilizes $B$.

We recall some basic terminology on root systems.
The {\it root system} of $(W,S)$ is the set $\Phi=\set{w(\a_s)\st w\in W,
s\in S}$ and the elements of
$\D$ are the {\it simple roots}.  The roots contained in
$$
 \Phi^+ = \big(\sum_{\a\in \D} \RM^+ \a\big) \cap \Phi
$$
are said to be {\it positive}, while those contained in
$\Phi^-=-\Phi^+$ are said to be {\it negative}. Moreover, $\Phi$
is the disjoint union of $\Phi^+$ and $\Phi^-$. If $w\in W$,
$\ell(w)=\card{N(w)}$, where
$$
N(w)=\set{\a \in \Phi^+ \st \ w(\a)\in \Phi^-}.
$$

Let $\a = w(\a_s)\in \Phi$, then $s_\a = wsw^{-1}$ acts as the
reflection in the hyperplane orthogonal to $\alpha$
and $s_\a = s_{-\a}$.
Therefore, the {\it set of reflections of $W$}
$$
T=\bigcup_{w\in W} wSw^{-1}
$$
is in bijection with $\Phi^+$ (and thus $\Phi^-$).

Let us recall the following well-known result:

\begin{lem}\label{important}
Let $w \in W$. Then:
\begin{itemize}
\itemth{a} If $\a \in \Phi^+$, then $\ell(ws_\a)> \ell(w)$ if and only if
$w(\a) \in \Phi^+$.
\itemth{b} If $s \in S$, then
$$N(sw)=\begin{cases}
        N(w) \coprod \{w^{-1}(\a_s)\} & \text{if } \ell(sw) > \ell(w), \\
    N(w) \setminus \{-w^{-1}(\a_s)\} & \text{otherwise.}
    \end{cases}$$
\end{itemize}
\end{lem}

From now on, and until the end of this paper, we fix a subset $A$
of $T$. We start with easy observations.

As a consequence of Lemma \ref{important} (a), we get that
$$
D_A(w)=\set{s_\a \in A \st \a \in \Phi^+ \text{ and } w(\a)\in
\Phi^-}.
$$
We also set
$$
N_A(w)=\{\a \in \Phi^+~|~s_\a \in A \text{ and } w(\a) \in
\Phi^-\}.
$$
The map $N_A(w) \to D_A(w)$, $\a \mapsto s_\a$ is then a
bijection.

\subsection{Properties of the map ${\boldsymbol{D_A}}$}
First, using Lemma \ref{important} (b), we get:

\begin{cor}\label{trivial}
If $s \in S$ and if $w \in W$ is such that $w^{-1} s w=s_{w^{-1}(\a_s)} \not\in A$,
then $N_A(w)=N_A(sw)$ (and $D_A(w)=D_A(sw)$).
\end{cor}

\remarque{A1 in A2} If  $A_1\subset A_2 \subset T$, then $D_{A_1} (w)=D_{A_2}(w)\cap A_1$
for all $w\in W$. Therefore if $W$ is finite, $\UD_{A_1}(W)
\subset \UD_{A_2}(W)$.~$\square$

\begin{prop}\label{A-admissible}
We have:
\begin{enumerate}
\item[(a)] $\varnothing$ is $A$-admissible.
\item[(b)] $D_\varnothing^A = \set{1}$ if and
only if  $S\subset A$.
\end{enumerate}
\end{prop}

\begin{proof}
We have $D_A(1)=\varnothing$ so (a) follows. If $s \in S\setminus
A$, then $D_A(s)=\varnothing$. This shows (c).
\end{proof}

The notion of $A$-descent set is obviously compatible with direct
products:

\begin{prop}\label{decomposable}
Assume that $W=W_1 \times W_2$ where $W_1$ and $W_2$ are standard
parabolic subgroups of $W$. Then, for all $I \in \PC_\ad(A)$, we
have
$$D_I^A = D_{I \cap W_1}^{A \cap W_1} \times D_{I \cap W_2}^{A \cap W_2}.$$
\end{prop}

\begin{cor}\label{produit algebres}
Assume that $W$ is finite and that $W=W_1 \times W_2$ where $W_1$
and $W_2$ are standard parabolic subgroups of $W$. Then
$$\UD_A(W) = \UD_{A \cap W_1}(W_1) \otimes_\ZM \UD_{A \cap W_2}(W_2).$$
\end{cor}

\exemple{group algebra} Consider the case where $A=T$ (then $N_A(w)=N(w)$).
It is well-known \cite[Chapter VI, Exercise 16]{bourbaki} that the map
$w \mapsto N(w)$ from $W$ onto the set of subsets of $\Phi^+$ is
injective  (observe that if $\alpha\in N(w_1w_2^{-1})$ then $\pm
w_2^{-1}(\alpha)$ lives in the union, but not in the intersection,
of $N(w_1)$ and $N(w_2)$). Therefore, the map $W
\to
\PC_\ad(T)$,
$w \mapsto D_T(w)$ is injective. In particular, if $W$ is finite,
then $\UD_T(W)=\ZM W$.\finl

\medskip

In the case of finite Coxeter groups, the multiplication on the
left by the longest element has the following easy property.

\begin{prop}\label{w0}
If $W$ is finite and if $w \in W$, then $D_A(w_0 w) = A \setminus
D_A(w)$.
\end{prop}

\begin{cor}\label{A-admissible-finite}
If $W$ is finite, then:
\begin{enumerate}
\item[(a)] $A$ is $A$-admissible;
\item[(b)]  $I\in \UP_\ad(A)$ if and only if $A\setminus I\in
\UP_\ad(A)$;
\item[(c)]  $D_A^A=\set{w_0}$ if and only
if $S\subset A$.
\end{enumerate}
\end{cor}

\begin{proof}
$D_A(w_0)=A$ so (a) follows. (b) follows from Proposition
\ref{w0}. (c) follows from Proposition \ref{w0} and Proposition
\ref{A-admissible} (b).
\end{proof}





%
%

\subsection{Left-connectedness}

In \cite{atkinson}, Atkinson gave a new proof of Solomon's result
by using an equivalence relation to describe descent sets. We
extend his result to $A$-descent sets. It shows in particular that
the subsets $D_I^A$ are left-connected.

Let $w$ and $w'$ be two elements of $W$. We say that $w$ is an
\textit{$A$-descent neighborhood} of $w'$, and we write $w
\smile_A w'$,
if $w' w^{-1} \in S$ and $w^{-1}w' \notin A$. It is easily seen
that $\smile_A$ is a symmetric relation. The reflexive and
transitive closure of the $A$-descent neighborhood relation is
called \textit{the $A$-descent equivalence}, and is denoted by
$\sim_A$. The next proposition characterizes this equivalence
relation in terms of $A$-descent sets.

\begin{prop}\label{descent_class}
Let $w$, $w' \in W$. Then
$$
w\edclass w' \iff D_A(w)=D_A(w') \iff N_A(w)=N_A(w').
$$
\end{prop}

\begin{proof}
The second equivalence is clear. If $w \smile_A w'$, then it
follows from Corollary \ref{trivial} that $N_A(w)=N_A(w')$. It
remains to show that, if $N_A(w)=N_A(w')$, then $w \sim_A w'$.

So, assume that $N_A(w)=N_A(w')$. Write $x=w' w^{-1}$ and let
$m=\ell(x)$. If $\ell(x)=0$, then $w=w'$ and we are done. Assume
that $m \ge 1$ and write $x=s_1 s_2 \dots s_m$ with $s_i \in S$.
We now want to prove on induction on $m$ that
$$w \smile_A s_m w \smile_A s_{m-1} s_m w \smile_A \dots \smile_A s_2 \dots s_m w
\smile_A s_1 s_2 \dots s_m w=w'.\leqno{(*)}$$
First, assume that $w \not\smile_A s_m w$. In other words, $w^{-1}
s_m w \in A$. For simplification, let $\a_i=\a_{s_i}$. By Lemma
\ref{important} (b), we have
$$N_A(s_m w) = N_A(w) \coprod \{w^{-1}(\a_m)\} \quad \text{or} \quad
N_A(w) \setminus \{-w^{-1}(\a_m)\}.
$$
In the first case, as $N_A(w)=N_A(w')$, and by applying again
Lemma \ref{important} (b), there exists a step $i \in
\{1,2,\dots,m-1\}$ between $N_A(w)$ to $N_A(w')$ where
$w^{-1}(\a_m)$ is removed from $N_A(s_i\dots s_m w)$, that is,
$(s_{i+1}\dots s_m w)^{-1}(\a_i)=-w^{-1}(\a_m)$. In the same way,
we get the same result in the second case. In other words, we have 
proved that there exists $i \in \{1,2,\dots,m-1\}$ such that 
$s_m\dots s_{i+1}(\a_i)=-\a_m$, so, by Lemma
\ref{important} (a), we have $\ell(s_m \dots s_{i+1} s_i) <
\ell(s_m \dots s_{i+1})$. This contradicts the fact that
$m=\ell(x)$. So $w \smile_A s_m w$, and  then $N_A(w)=N_A(s_m w)$.
Hence $N_A(s_m w)=N_A(w')$ and $w' (s_m w)^{-1} = s_1 s_2 \dots
s_{m-1}$.  We get then by induction that
$s_m w \smile_A s_{m-1} s_m w \smile_A \dots \smile_A s_1 s_2
\dots s_m w=w'$, which shows $(*)$.
\end{proof}

\begin{cor}
If $I \in \PC_\ad(A)$, then $D_I^A$ is left-connected.
\end{cor}

\subsection{Nice subsets of ${\boldsymbol{T}}$}

We say that $A$ is {\it nice} if, for every $s \in A$ and $w \in
W$ such that $w^{-1} s w \not\in A$, we have $D_A(sw)=D_A(w)$.
Notice that every subset of $S$ is nice, by
Corollary~\ref{trivial}.

If $w\in W$ and $I$, $J\in\UP_\ad (A)$, we set
$$
D_A(I,J,w) = \set{(u,v) \in D_I^A \times D_J^A\st uv=w }.
$$
The next lemma gives an obvious characterization of the fact
that $\DC_A(W)$ is an algebra in terms of these sets.

\begin{lem}\label{key}
Assume that $W$ is finite. Then the following are equivalent:
\begin{itemize}
\itemth{1} $\UD_A(W)$ is a subalgebra of $\ZM W$.
\itemth{2} For all $I$, $J \in \PC_\ad(A)$ and for all $w$, $w' \in W$
such that $D_A(w)=D_A(w')$, we have $|D_A(I,J,w)|=|D_A(I,J,w')|$.
\end{itemize}
If these conditions are fulfilled,
we choose for any $I\in\PC_\ad (A)$ an element $z_I$ in $D_I^A$.
Then
$$d_I^A d_J^A = \sum_{K \in \PC_\ad(A)} |D_A(I,J,z_K)| d_K^A.$$
\end{lem}

\def\Id{{\mathrm{Id}}}

Let us fix now $s \in S$ and  let $(u,v) \in W \times W$. If $u
\smile_A su$, we set $\psi_s^A(u,v)=(su,v)$. If
$u \not\smile_A su$, then we set $\psi_s^A(u,v)=(u,u^{-1}su v)$.
Note that, in the last case, $u^{-1}su \in A$. We have
$(\psi_s^A)^2 = \Id_{W \times W}$. In particular, $\psi_s^A$ is a
bijection. Using $\psi_s^A$, one can relate the notion of nice
subsets to the property (2) stated in Lemma \ref{key}.

\begin{prop}\label{key 2}
Assume that $A$ is nice. Let $I$, $J \in \PC_\ad(W)$, let $w \in
W$ and let $s \in S$ be such that $w \smile_A s w$. Then
$\psi_s^A(D_A(I,J,w))=D_A(I,J,sw)$.
\end{prop}

\begin{proof}
Let $(u,v)$ be an element of $D_A(I,J,w)$. By symmetry, we only need to prove
that $\psi_s^A(u,v) \in D_A(I,J,sw)$. If $u \smile_A su$, then
$D_A(su)=D_A(u)=I$ by Proposition \ref{descent_class}, so
$\psi_s^A(u,v)=(su,v) \in D_A(I,J,sw)$. So, we may, and we will,
assume that $u \not\smile_A su$. Let $s'=u^{-1} s u \in A$. Note
that $w^{-1} s w =v^{-1}s'v \not\in A$. Then,
$\psi_s^A(u,v)=(u,s'v)$ and $us'v=sw$. So, we only need to prove
that $D_A(s'v)=D_A(v)$. But this just follows from the definition
of nice subset of $T$.
\end{proof}

\begin{cor}\label{voila}
If $W$ is finite and if $A$ is nice, then $\UD_A(W)$ is a
subalgebra of $\ZM W$. It is unitary if and only if $S \subset A$.
\end{cor}

\begin{proof}
This follows from Lemma \ref{key} and from Propositions
\ref{descent_class} and \ref{key 2}.
\end{proof}

\remarque{question nice} It would be interesting to know if the
converse of Corollary \ref{voila} is true.\finl

\subsection{Proof of Theorems A and B}

Using Corollary \ref{voila}, we see that Theorems A and B are
direct consequences of the following theorem (which holds also for
infinite Coxeter groups):

\begin{thm}\label{le theoreme}
Assume that one of the following holds:
\begin{itemize}
\itemth{1} There exists two subsets $S_1$ and $S_2$ of $S$ such that
$A=S_1 \cup (\cup_{s \in S_2} C(s))$.
\itemth{2} $S=\{s,t\}$, $m(s,t)$ is even or $\infty$,
and $A=\set{s,t,sts}$ or $A=\set{t,sts}$.
\end{itemize}
Then $A$ is nice.
\end{thm}

\begin{proof}
Assume that (1) or (2) holds. Let $r \in A$ and let $w \in W$ be
such that $w^{-1} r w \not\in A$. We want to prove that
$D_A(rw)=D_A(w)$. By symmetry, we only need to show that $D_A(rw)
\subset D_A(w)$. If $r \in S$, then this follows from Corollary
\ref{trivial}. So we may, and we will, assume that $r \notin S$.

$\bullet$ Assume that (1) holds. Write $A'=\cup_{s \in S_2} C(s)$
and $S' = A \setminus A'$. Then $S' \subset S$, $A=A' \coprod S'$
and $A'$ is stable under conjugacy. Then $r\in A'$ and $w^{-1} r w
\in A'\subset A$, which contradicts our hypothesis.

$\bullet$ Assume that (2) holds. If $m(s,t)=2$, then $A$ is
contained in $S$ and therefore is nice by Corollary~\ref{trivial}.
So we may, and we will, assume that
$m(s,t)\ge 4$. Since $r \not\in S$, we have $r=sts$.
Assume that $D_A(rw)\not\subset D_A(w)$. Let
$$\Phi_A=\{\a \in \Phi^+~|~s_\a \in A\} \subset \{\a_s,\a_t,s(\a_t)\}.$$
There exists $\a \in \Phi_A$ such that $rw(\a) \in \Phi^-$ and
$w(\a) \in
\Phi^+$. So $w(\a) \in N(r)=\{\a_s,s(\a_t),st(\a_s)\}$. Since
$w^{-1} r w \not\in A$, we have that $w s_\a w^{-1} \not= r$, so
$w(\a) \not= s(\a_t)$. So $w(\a)=\a_s$ or $st(\a_s)$. But the
roots $\a_s$ and $\a_t$ lie in different $W$-orbits. So $\a=\a_s$
and $w(\a_s) \in \{\a_s,st(\a_s)\}$. If $W$ is infinite, this
gives that $w \in \{1,st\}$. This contradicts the fact that
$w^{-1} rw \not\in A$. If $W$ is finite, this gives that $w \in
\{1,st,w_0s,sts w_0\}$. But again, this contradicts the fact that
$w^{-1} rw \not\in A$.
\end{proof}

\exemple{f4}
Assume here that $W$ is of type $F_4$, that
$S=\{s_1,s_2,s_3,s_4\}$ and that $A=C(s_1) \cup S$. Then, using
{\tt GAP}, one can see that $\rank_\ZM \UD_A(W)=300$ and
$\rank_\ZM \Sigma_A(W)=149$. Moreover, $\Sigma_A(W)$ is not a
subalgebra of $\UD_A(W)$.\finl

\medskip

Some computations with {\tt GAP} suggest that the following
question has a positive answer:

\smallskip

\begin{quotation}
\noindent{\bf Question. } {\it Let $A$ be a nice subset of $T$
containing $S$. Is it true that $A$ is one of the subsets
mentioned in Theorem \ref{le theoreme}?}
\end{quotation}


%
%

\def\sba{{\bar{s}}}
\def\tba{{\bar{t}}}

\section{The dihedral groups}\label{section 2}

The aim of this section is to study the unitary subalgebras
$\UD_A(W)$ constructed in Theorems A and B whenever $W$ is finite
and dihedral.

\smallskip

\begin{quotation}
\noindent{\bf Hypothesis:} {\it From now on, and until the end of
this paper, we assume that $S=\{s,t\}$ with $s \not= t$ and that
$m(s,t)=2m$, with $2 \le m < \infty$.}
\end{quotation}

\smallskip

Note that $w_0=(st)^m$ is central. We denote by $\PC_0(A)$ the set
of subsets $I$ of $A$ such that $W_I \cap A = I$. With this
notation, we have
$$\Sigma_A(W)=\sum_{I \in \PC_0(A)} \ZM x_I.$$
We now define an equivalence relation $\equiv$ on $\PC_0(A)$: if
$I$, $J \in \PC_0(A)$, we write $I \equiv J$ if $W_I$ and $W_J$
are conjugate in $W$. We set
$$\Sigma^{(1)}_A(W)=
\sum_{\SS{I,J \in \PC_0(A)} \atop \SS{I \equiv J}} \ZM(x_I-x_J).$$

In what follows, we will also need some facts on the character
table of $W$. Let us recall here the construction of $\Irr W$.
First, let $H$ be the subgroup of $W$ generated by $st$. It is
normal in $W$, of order $2m$ (in other words, of index $2$). We
choose the primitive $(2m)$-th root of unity $\zeta \in \CM$ of
argument $\frac{\pi}{m}$. If $i \in
\ZM$, we denote by $\xi_i: H \to \CM^\times$ the unique linear
character such that $\xi_i(st)=\z^i$. Then $\Irr H = \{\xi_i~|~0
\le i \le 2m-1\}$. Now, let
$$
\chi_i=\Ind_H^W \xi_i.
$$
Then $\chi_i=\chi_{2m-i}$ and, if $1 \le i \le m-1$, $\chi_i \in
\Irr W$. Also, $\chi_i$ has values in $\RM$. More precisely, for
$1 \le i \le m-1$ and $j\in\ZM$
$$
\chi_i\big((ts)^j\big)=\z^j + \z^{-j}=2\cos\left(\frac{ij\pi}{m}\right)
\quad\textrm{and}\quad\chi_i\big(s(ts)^j\big)=0.
$$
Let $1$ denote the trivial character of $W$, let $\e$ denote the
sign character and let $\g: W \to \{1,-1\}$ be the unique linear
character such that $\g(s)=-\g(t)=1$. Then
\equat\label{irreductibles W} \Irr W = \{1,\e,\g,\e\g\} \cup
\{\chi_i~|~1 \le i \le m-1\}.
\endequat
In particular, $|\Irr W|=m+3$.

\subsection{The subset ${\boldsymbol{A=\{s,t,sts\}}}$}\label{A}

From now on, we assume that $A=\{s,t,sts\}$.  We set
$\sba=A \setminus \{s\}=\{t,sts\}$ and $\tba=A
\setminus\{t\}=\{s,sts\}$. It is easy to see that
$$\PC_\ad(A)=\{\varnothing,\{s\},\{t\},\sba,\tba,A\}.$$
For simplification,  we will denote by $d_I$ the element $d_I^A$ of
$\ZM W$ (for $I \in \PC_\ad(W)$) and we set $d_s=d_{\{s\}}$ and
$d_t=d_{\{t\}}$. We have
$$
\begin{array}{rcl@{\qquad}rcl}
d_\varnothing &=& 1, & d_\sba&=& w_0 s, \\
d_s   &=& s, &d_A   &=& w_0, \\
d_t   &=&\DS{\sum_{i=1}^{m-1}} \Bigl((st)^i + (ts)^{i-1}t\Bigr), &
d_\tba&=&\DS{\sum_{i=1}^{m-1}} \Bigl((st)^is + (ts)^i\Bigr). \\
\end{array}
$$
The multiplication table of $\UD_A(W)$ is given by

\medskip
{\small
\begin{centerline}
{\begin{tabular}{|c||c|c|c|c|c|c|}
\hline
\vertical         & $1$      & $d_s$    & $d_\sba$ & $d_A$    & $d_t$    & $d_\tba$ \\
\hline
\hline
\vertical $1$      & $1$      & $d_s$    & $d_\sba$ & $d_A$    & $d_t$    & $d_\tba$ \\
\vertical $d_s$    & $d_s$    & $1$      & $d_A$    & $d_\sba$ & $d_t$    & $d_\tba$ \\
\vertical $d_\sba$ & $d_\sba$ & $d_A$    & $1$      & $d_s$    & $d_\tba$ & $d_t$    \\
\vertical $d_A$    & $d_A$    & $d_\sba$ & $d_s$    & $1$      & $d_\tba$ & $d_t$    \\
\vertical $d_t$    & $d_t$    & $d_\tba$ & $d_t$    & $d_\tba$ & $z_A$    & $z_A$    \\
\vertical $d_\tba$ & $d_\tba$ & $d_t$   & $d_\tba$ & $d_t$    & $z_A$    & $z_A$    \\
\hline
\end{tabular}}
\end{centerline}
}

\medskip

\noindent where
$z_A=(m-1)\big(1+d_A+d_s+d_\sba\big)+(m-2)\big(d_t+d_\tba\big)$.
We now study the sub-$\ZM$-module $\Sigma_A(W)$: we will show that
it coincides with $\UD_A(W)$. First, it is easily seen that
$$
\UP_0 (A)=\set{\varnothing,\, \set{s},\, \set{t},\, \set{sts},\, \sba,\, A}
$$
and that
$$\begin{array}{ccclllll}
x_A    & = & 1 &&&&& \\
x_\sba & = & 1 & + d_s &&&& \\
x_{sts} & = & 1 & + d_s & + d_t &&& \\
x_t    & = & 1 & + d_s &       & +d_\tba && \\
x_s    & = & 1 &       & + d_t &         & + d_\sba & \\
x_\varnothing&=&1&+d_s & + d_t & +d_\tba & + d_\sba & + d_A \\
\end{array}$$
Therefore, $\Sigma_A(W)=\UD_A(W)=\oplus_{I \in \PC_0(A)} \ZM x_I$.
So we can define a map $\th_A: \Sigma_A(W) \to \ZM \Irr W$ by
$\th_A(x_I)=\Ind_{W_I}^W 1_{W_I}$.

\begin{prop}\label{s,t,sts}
Assume that $S=\{s,t\}$ with $s \not= t$, that $m(s,t)=2m$ with $m
\ge 2$ and that $A=\{s,t,sts\}$. Then:
\begin{itemize}
\itemth{a} $\Sigma_A(W)=\UD_A(W)$ is a subalgebra of $\ZM W$ of $\ZM$-rank $6$.
\itemth{b} $\th_A$ is a morphism of algebras.
\itemth{c} $\Ker \th_A = \ZM (x_t - x_{sts}) = \Sigma_A^{(1)}(W)$.
\itemth{d} $\QM \otimes_\ZM \Ker \th_A$ is the
radical of $\QM \otimes_\ZM \Sigma_A(W)$.
\itemth{e} $\th_A$ is surjective if and only if $m=2$ that is, if and only if
$W$ is of type $B_2$.
\end{itemize}
\end{prop}

\begin{proof}
(a) has already been proved.
For proving the other assertions, we need to compute explicitly
the map $\th_A$. It is given by the following table:
$$\begin{array}{|c||c|c|c|c|c|c|}
\hline
\verticala d_I & 1 & d_s & d_\sba & d_A & d_t & d_\tba \\
\hline
\verticala \th_A(d_I)&1&\e\g&\g&\e&\DS{\sum_{i=1}^{m-1}}\chi_i&\DS{\sum_{i=1}^{m-1}}\chi_i\\
\hline
\end{array}$$

(c) This shows that $\Ker \th_A
=\ZM(d_\tba-d_t)=\ZM(x_t-x_{sts})=\Sigma_A^{(1)}(W)$, so (c)
holds.

(b) To prove that $\th_A$ is a morphism of algebras, the only
difficult point is to prove that $\th_A(d_t^2)=\th_A(d_t)^2$. Let
$\r$ denote the regular character of $W$. Then
$$\th_A(d_t)=\frac{1}{2}\big(\r-1-\g-\e\g-\e).$$
Therefore,
$$\th_A(d_t)^2= (m-2) \r + 1+\g+\e\g+\e.$$
But, $d_t^2=z_A=(m-2) x_\varnothing + 1 + d_s + d_\sba + d_A$.
This shows that $\th_A(d_t^2)=\th_A(d_t)^2$.

(d) Let $R$ denote the radical of $\QM \otimes_\ZM \Sigma_A(W)$.
We only need to prove that $\CM
\otimes_\QM R = \CM \otimes_\ZM \Sigma_A(W)$. Since $\CM\Irr W$ is
a split semisimple commutative algebra, every subalgebra of
$\CM\Irr W$ is semisimple. So $(\CM
\otimes_\ZM \Sigma_A(W))/(\CM \otimes_\ZM \Ker \th_A)$ is a
semisimple algebra. This shows that $R$ is contained in $\CM
\otimes_\ZM \Ker \th_A$. Moreover, since
$(x_t-s_{sts})^2=(d_t-d_\tba)^2=0$, $\Ker \th_A$ is a nilpotent
two-sided ideal of $\Sigma_A(W)$. So $\CM \otimes_\ZM \Ker \th_A$
is contained in $\CM \otimes_\ZM R$. This shows (d).

(e) If $m=2$, then $\Irr W=\{1,\g,\e\g,\e,\chi_1\}
=\th_A(\{1,d_s,d_\sba,d_A,d_t\})$ so $\th_A$ is surjective.
Conversely, if $\th_A$ is surjective, then $|\Irr W|=5$ (by (a)
and (c)). Since $|\Irr W|=m+3$, this gives $m=2$.
\end{proof}

We close this subsection by giving, for $\Sigma_A(W)$, a complete set
of orthogonal primitive idempotents,  extending to our case those
given in \cite{bbht}, and the corresponding  irreducible
representations:
%
%
%
%
$$
E_\varnothing=\frac{1}{4m}x_\varnothing, \qquad
E_s=\frac{1}{2}\left(x_s-\frac{1}{2}x_\varnothing\right), \qquad
E_t=\frac{1}{2}\left(x_t-\frac{1}{2}x_\varnothing\right),
$$
\begin{eqnarray*}
E_{\bar s}&=&\frac{1}{2}\left(x_{\bar s}-\frac{1}{2}x_t-\frac{1}{2}x_{sts}+\frac{m-1}{2m}x_\varnothing\right)\\
E_A&=&1-\frac{1}{2}x_s-\frac{1}{4}x_t+\frac{1}{4}x_{sts}-\frac{1}{2}x_\sba+\frac{1}{4}x_\varnothing
\end{eqnarray*}

\def\ev{{\mathrm{ev}}}

For $i \in \mathcal  \{\varnothing, s,t,\sba,A\}$, we denote by
$P_i$ the indecomposable projective module $\QM\Sigma_A(W) E_i$.
Let $\QM_i^A$ denote the unique simple $\QM\Sigma_A(W)$-module
lying in the head of $P_i$. It is easily seen, using the previous
multiplication table, that
$$\dim P_\varnothing=\dim P_t = \dim P_\sba = \dim P_A=1$$
$$\dim P_s =2\leqno{\text{and that}}$$
(note that $x_\sba E_s = -(x_t-x_{sts})$). In other words,
$$P_\varnothing=\QM_\varnothing^A,\quad P_t = \QM_t^A, \quad P_\sba=\QM_\sba^A,\quad
P_A=\QM_A^A$$
$$\Rad(P_s)\simeq \QM_t^A.\leqno{\text{and}}$$

For $w \in W$, we denote by $\ev_w$ the morphism of algebras
$\ZM\Irr W \to \ZM$, $\chi \mapsto
\chi(w)$ and we set $\ev_w^A = \ev_w \circ \th_A$. Then
the morphism of algebras
$\Sigma_A(W) \to \ZM$ associated to the
simple module $\QM_i^A$ is $\ev_{f(i)}^A$, where
$$f(\varnothing)=1,\quad f(s)=s,\quad f(t)=t, \quad f(\sba)= (st)^m = w_0
\quad{\text{and}}\quad f(A)=st.$$

\subsection{The subset $\boldsymbol{B=\set{s}\cup C(t)}$}

Let $B=\set{s}\cup C(t)$ (so that $|B|=m+1$).
It is easy to see that
$\PC_\ad(B)$ consists of the sets
$$
\begin{array}{ccl}
\varnothing                &=& D_B(1), \\
B                        &=& D_B(w_0), \\
\set{s}                  &=& D_B(s), \\
C(t)                     &=& D_B (w_0s), \\
D_B\big((ts)^i \big)     &=& D_B\big( s(t s)^i\big), \quad 1\leq i\leq m-1,\\
D_B\big((st)^{j}\big) &=& D_B\big( (t s)^{j-1} t\big), \quad 1\leq j\leq m-1. \\
\end{array}
$$
Therefore $\UD_B(W)$ is a subalgebra of $\ZM W$ of $\ZM$-rank
$(2m+2)$.

Using {\tt GAP}, we can see that, in general, $\Sigma_B
(W)\not=\UD_B(W)$. First examples are given in the following table
$$\begin{array}{|c||c|c|c|c|c|c|c|c|c|c|c|}
\hline
\verticala  m                                   & 2 & 3 & 4 & 5 & 6 & 7 & 8 & 9 & 10 & 11    \\
\hline
\verticala \ZM\textnormal{-rank of }\UD_B(W)    & 6 & 8 & 10& 12& 14& 16& 18& 20& 22 & 24    \\
\hline
\verticala \ZM\textnormal{-rank of }\Sigma_B(W) & 6 & 8 & 10& 10& 14& 12& 18& 18& 22 & 16 \\
\hline
\end{array}$$

\remarque{surjectif}
The linear map $\theta_B:\Sigma_B(W)\to \ZM\Irr
W$, $x_I\mapsto \Ind_{W_I}^W 1_{W_I}$, is well-defined and
surjective if and only if $m \in \{2,3\}$ (recall that $m \ge 2$).
Indeed, the image of $\theta_B$ can not contain a non-rational
character. But all characters of $W$ are rational if and only if $W$
is a Weyl group, that is, if and only if $2m\in\set{2,3,4,6}$.

Moreover, if $m=2$, then $A=B$ and $\th_B=\th_A$ is surjective by
Proposition~\ref{s,t,sts}. If $m=3$, then this follows from
Proposition \ref{s,t,sts,tstst} below.\finl

\subsection{The algebra $\boldsymbol{\UD_B(G_2)}$}

From now on, and until the end, we assume that $m=3$. That is $W$
is of type $I_2(6)=G_2$. For convenience, we keep the same
notation as in \S\ref{A}. We have
$$\begin{array}{rcl}
d_\varnothing &=& 1,\\
d_s         &=& s,\\
d_\sba      &=& d^B_{C(t)}            = w_0 s,\\
d_A         &=& d^B_B                 = w_0,
\end{array}
\qquad
\begin{array}{rcl}
d_1         &=& d^B_\set{t}           =t+st,\\
d_2         &=& d^B_\set{s,sts}       = ts+sts,\\
d_3         &=& d^B_\set{s,sts,tstst} = tsts+ststs,\\
d_4         &=& d^B_\set{t,tstst}     = tst+stst.
\end{array}
$$

Let us now show that $\Sigma_B(W)=\UD_B(W)$. First, it is easily seen that
$$
\UP_0 (B)=\set{\varnothing,\, \set{s},\, \set{t},\, \set{sts},\,
\sba,\, \set{tstst},\set{s,tstst},B}
$$
and that
$$\begin{array}{ccclllllll}
x_A = x_B            & = & 1 &&&&& &&\\
x_\sba              & = & 1 & + d_s &&&&&& \\
x_\set{s,tstst} & = & 1 &       &+d_1 &&&&&\\
x_{tstst} & = & 1 &  + d_s&+d_1 & + d_2  &&&&     \\
x_t                 & = & 1 & +d_s  &     & +d_2       & +d_3&&&\\
x_{sts}             & = & 1 & + d_s & + d_1 &&&+d_4&& \\
x_s                 & = & 1 &       & + d_1 &         &  &+d_4&+d_\sba& \\
x_\varnothing       & = &1&+d_s & + d_1 &+d_2& +d_3 & + d_4 & + d_\sba &+d_A\\
\end{array}$$
Therefore, $\Sigma_B(W)=\UD_B(W)=\oplus_{I \in \PC_0(B)} \ZM
x_I$. So we can define a map $\th_B: \Sigma_B(W) \to \ZM \Irr W$
by $\th_B(x_I)=\Ind_{W_I}^W 1_{W_I}$.


\begin{prop}\label{s,t,sts,tstst}
Assume that $S=\{s,t\}$ with $s \not= t$, that $m(s,t)=6$ and that
$B=\{s,t,sts,tstst\}$. Then:
\begin{itemize}
\itemth{a} $\Sigma_B(W)=\UD_B(W)$ is a subalgebra of $\ZM W$ of
$\ZM$-rank $8$.
\itemth{b} $\th_B$ is a surjective  linear map (and not a morphism
of algebras).
\itemth{c} $\Ker \th_B = \ZM (x_{tstst} -x_t)\oplus\ZM(x_{tstst}-
x_{sts}) = \Sigma_B^{(1)}(W)$.
\itemth{d} $\Irr W=\theta_B\left(\set{1,d_s,d_\sba,d_A,d_1,d_2}\right)$.
\end{itemize}
\end{prop}

\begin{proof}
The map $\th_B$ is given by the following table:
$$\begin{array}{|c||c|c|c|c|c|c|c|c|}
\hline
\verticala d^B_I & 1 & d_s & d_\sba & d_A & d_1 & d_2 & d_3&d_4 \\
\hline
\verticala \th_B(d^B_I)&1&\e\g&\g&\e&\chi_2&\chi_1&\chi_2 &\chi_1 \\
\hline
\end{array}$$

This shows (b) and (d). As $d_1-d_3=x_\set{tstst}-x_t$ and
$d_2-d_4=x_{tstst}-x_{sts}$, (c) is proved. Finally, the fact that
$\theta_B$ is not a morphism of algebras follows from
$\th_B(d_1^2)(w_0)=(1+\e\g+\chi_1)(w_0)=-2\not=
\th_B(d_1)^2(w_0)=\chi_2(w_0)^2=4.$
\end{proof}

For information,
the multiplication table in the basis $(d_I^B)_{I \in \UP_\ad(B)}$ is given by

\medskip

{\small
\begin{centerline}
{\begin{tabular}{c}
\begin{tabular}{|c||c|c|c|c|c|c|c|c|}
\hline
\vertical             & $1$        & $d_s$      & $d_1$         & $d_2$         & $d_3$                  & $d_4$                   & $d_{\sba}$ & $d_A$      \\
\hline
\hline
\vertical  $1$        & $1$        & $d_s$      & $d_{1}$       & $d_2$         & $d_3$                  & $d_4$                   & $d_{\sba}$ & $d_A$      \\
\vertical  $d_s$      & $d_s$      & $1$        & $d_{1}$       & $d_2$         & $d_3$                  & $d_4$                   & $d_A$      & $d_{\sba}$ \\
\vertical  $d_1$      & $d_{1}$    & $d_2$      & $1 + d_s + d_4$ & $1 + d_s + d_3$ & $d_2 + d_{\sba} + d_A$ & $d_1 + d_{\sba} + d_A$  & $d_4$  & $d_3$      \\
\vertical  $d_2$      & $d_2$      & $d_{1}$    & $1 + d_s + d_4$ & $1 + d_s + d_3$ & $d_2 + d_{\sba} + d_A$ & $d_1 + d_{\sba} + d_A$  & $d_3$  & $d_4$      \\
\vertical  $d_3$      & $d_3$      & $d_4$      & $d_2 + d_{\sba} + d_A$ & $d_1 + d_{\sba} + d_A$ & $1 + d_s + d_4$ & $1 + d_s + d_3$  &  $d_2$   & $d_{1}$    \\
\vertical  $d_4$      & $d_4$      & $d_3$      & $d_2 + d_{\sba} + d_A$ & $d_1 + d_{\sba} + d_A$ & $1 + d_s + d_4$ & $1 + d_s + d_3$  & $d_1$    & $d_2$      \\
\vertical  $d_{\sba}$ & $d_{\sba}$ & $d_A$      & $d_3$         & $d_4$         &  $d_{1}$               & $d_2$                   & $1$      &  $d_s$       \\
\vertical  $d_A$      & $d_A$      & $d_{\sba}$ & $d_3$         & $d_4$         &  $d_{1}$               & $d_2$                   & $d_s$        &  $1$     \\
\hline
\end{tabular}\\
\end{tabular}}
\end{centerline}
}

\end{document}